\documentclass[12pt,leqno]{article}  
\usepackage{a4wide,epsf,latexsym}
\usepackage[centertags]{amsmath}
\usepackage{graphicx}
\usepackage{cite}
\numberwithin{equation}{section}

\newtheorem{remark}{Remark}

\begin{document}

\title {Robust Numerical Methods for Singularly Perturbed Differential
Equations--Supplements
}

\author{Hans-G. Roos, TU Dresden}

\date{2022}

\maketitle

\begin{abstract}
The second edition of the book \cite{RST08} appeared many years ago and was for many years a reliable guide into the world of
numerical methods for singularly perturbed problems. Since then many new results came into the game,
we present some selected ones and the related sources.
\end{abstract}

{\it AMS subject classification}: 65 L10, 65 L12, 65 L50, 65 N30, 65 N50

\section{Introduction}
Since the appearance of \cite{RST08} many new results improved our understanding of the
complicated task to solve singularly perturbed and flow problems numerically.
The survey \cite{Ro12a} discussed new
developments in the years 2008-2012,  and in a paper 2015 \cite{RS15a} we discussed some open
problems in analysis and numerics of singularly perturbed problems. Recently,  John,  Knobloch
and  Novo also published a  survey on recent numerical methods for convection-dominated
equations and incompressible flow problems and sketched some open problems \cite{JKN18}.

Here we present some selected topics which are not our not so detailed to find in \cite{RST08} and
hints to related publications.

\section{Systems of ordinary differential equations}
Systems of ordinary differential equations are often discussed in books on asymptotic expansions
for singularly perturbed problems: see, e.g.,  \cite[Chapter 2]{VB90} or
\cite[Chapter 7]{Wa65}. Nevertheless in the past relatively little attention was paid to their
numerical solution, although the papers \cite{Ba69} (on reaction-diffusion
systems) and \cite{AKK74} (on convection-diffusion systems) are
worth noting. In recent years interest in this area has grown, as we now describe.

Consider a general system of $M$ equations:
\begin{subequations}\label{convreacdiff}
\begin{align}
L\vec{u}:&= -\varepsilon\vec{u}''+B\vec{u}' + A\vec{u}= \vec{f}
\quad\text{on }\Omega:=(0,1),\\
\vec{u}(0) &=\vec{g}_0,\ \vec{u}(1) =\vec{g}_1,
\end{align}
\end{subequations}
where $\vec{u}=(u_1,u_2,\dots,u_M)^T$ is the unknown solution
while $\vec{f} = (f_1,\dots, f_M)^T$, $\vec{g}_0$ and $\vec{g}_1$
are constant column vectors, and $A=(a_{ij})$ and $B=(b_{ij})$ are
$M\times M$ matrices.

The system \eqref{convreacdiff} is said to be \emph{weakly
coupled} if the convection coupling matrix $B$ is diagonal, i.e., the $i^{\rm th}$ equation of the
system is
\begin{equation}\label{convdiff-itheq}
 -\varepsilon u_i'' + b_{ii}u_i' + \sum_{j=1}^M a_{ij}u_j = f_i,
\end{equation}
so the system is coupled only through the lower-order reaction
terms.

Lin{\ss} \cite{Lin07a} allows different diffusion coefficients in
different equations: $\varepsilon = \varepsilon_i$ in the $i^{\rm
th}$ equation for $i=1,\dots, M$. Assume that $b_{ii}(x) \ge
\beta_i >0$ and $a_{ii}(x) \ge \alpha >0$ on $[0,1]$ for each $i$.
(In \cite{Lin07a} the weaker hypothesis $|b_{ii}(x)| \ge \beta_i
>0$ is used, which permits layers in $\vec{u}$ at both ends of [0,1], but for
brevity we won't consider this here.) Rewrite
\eqref{convdiff-itheq} as
\begin{equation}\label{convdiff-itheqB}
 -\varepsilon_i u_i'' + b_{ii}u_i' + a_{ii} u_i = - \sum_{j\ne i} a_{ij}u_j + f_i,
\end{equation}
Then  $ \|u_i\|_\infty \le \| (- \sum_{j\ne i} a_{ij}u_j +
f_i)/a_{ii} \|_\infty$ by a standard maximum principle argument.
Rearranging, one gets
$$ \|u_i\|_\infty - \sum_{j\ne i} \left\|\frac{a_{ij}}{a_{ii}}\right\|_\infty
        \|u_j\|_\infty \le  \left\|\frac{f_i}{a_{ii}}\right\|_\infty
    \quad\text{for } i=1,\dots, M. $$
Define the $M \times M$ matrix $\Gamma = (\gamma_{ij})$ by $\gamma_{ii} =1$, $\gamma_{ij} =
-\|a_{ij}/a_{ii}\|_\infty$ for $i \ne j$. Assume that $\Gamma$ is inverse-monotone, i.e., that
$\Gamma^{-1} \ge 0$. It follows that $\|\vec u\|_\infty \le C \|\vec f\|_\infty$ for some constant
$C$, where $\|\vec v\|_\infty = \max_i \|v_i\|_\infty$ for $\vec v = (v_1,\dots,v_M)^T$.  One can
now apply a scalar-equation analysis  for each
$i$ and get
\begin{equation}\label{syst,deriv}
 |u_i^{(k)}(x)| \le C\left[ 1 + \varepsilon_i^{-k}
    e^{-\beta_i (1-x)/\varepsilon_i} \right] \quad\text{for }x\in [0,1] \text{ and } k=0,1.
    \end{equation}

Thus there is no strong interaction between the layers in the
first-order derivatives of different components $u_i$;
nevertheless the domains of these layers can overlap and this
influences the construction of numerical methods for
\eqref{convreacdiff}.
But for higher order derivatives additionally weak layers come into the play. In \cite{RS15} it is
proved
$$ |u_i''(x)| \le C\left( 1 + \varepsilon_i^{-2}
    e^{-\beta(1-x)/\varepsilon_i}+\frac{1}{\delta} \sum_{j=i+1}^M\varepsilon_j^{-1}
    e^{-\beta(1-\delta
    )(1-x)/\varepsilon_j}\right),
    $$
here $\delta>0$ is arbitrary, assuming $b_{kk}\ge \beta>0$. The estimate shows that the $i$-th component
is overlapped by weak layers generated by other components, but only these for $j>i$ due to the
increase in diffusion $\varepsilon_1<\varepsilon_2<\cdots$.

The system \eqref{convreacdiff} is said to be \emph{strongly
coupled} if for some $i\in\{1,\dots,M\}$ one has $b_{ij}\ne 0$ for some $j\ne i$. Such systems do
not satisfy a maximum principle of the usual type. One now gets stronger interactions between
layers even in the case $\varepsilon_i=\varepsilon$
;
 see \cite{AKK74,Lin07b,ORS08,ORSS08,RR11}.

 We start with that case and the example
\begin{subequations}
\begin{align}
 -\varepsilon u_1''-3u_1'-4u_2' &=1 ,\nonumber\\
     -\varepsilon u_2''-4u_1'+3u_2'&=2\nonumber
\end{align}
\end{subequations}
and homogeneous boundary conditions at $x=0$ and $x=1$. Canceling of exponentially small terms yields
the asymptotic approximation
\begin{subequations}
\begin{align}
 u_1^{as} &=8/25-11/25x-8/25e^{-5x/\varepsilon}+3/25e^{-5(1-x)/\varepsilon},
     \nonumber\\
    u_2^{as}  &=4/25+2/25x-4/25e^{-5x/\varepsilon}-6/25e^{-5(1-x)/\varepsilon}.
    \nonumber
\end{align}
\end{subequations}
Both solution components do have strong layers at $x=0$ and $x=1$ which correspond to the
eigenvalues of the matrix $B$. Moreover, the reduced solution
\begin{subequations}
\begin{align}
 u_1^{0} &=8/25-11/25x,
     \nonumber\\
    u_2^{0}  &=4/25+2/25x
    \nonumber
\end{align}
\end{subequations}
does not satisfy any of the given boundary conditions, but it holds
$$ (u_1-2u_2)(0)=0,\quad (2u_1+u_2)(1)=0,
$$
which corresponds to the eigenvectors of $B$.

If we assume that all eigenvalues of $B$ do not change the sign in the given interval,
every component $u_i$ has, in general, $M$ layers but their location depend on the sign
pattern of the eigenvalues (assuming for simplicity, that $B$ is symmetric).

In the case that all eigenvalues have one sign, all overlapping boundary layers are located
at one boundary. Such a case is studied in \cite{ORS08,ORSS08} without assuming symmetry of
$B$. We follow \cite{ORSS08} to derive a stability result.

 For each $i$ assume $b_{ii}(x) \ge \beta_i>0$ and
$a_{ii}(x)\ge 0$ on $[0,1]$. Rewrite the $i^{\rm th}$ equation as
\begin{subequations}\label{convdiffstrong}
\begin{align}
L_i u := -\varepsilon u_i'' + b_{ii}u_i' + a_{ii}u_i
    &= f_i + \sum_{\genfrac{}{}{0pt}{}{j=1}{j\ne i}}^m \left[
    (b_{ij}u_j)' - (b_{ij}'+a_{ij})u_j \right],\\
    u_i(0) &= u_i(1) = 0.
\end{align}
\end{subequations}
For the scalar problem $L_iv= \phi$ and $v(0)=v(1)=0$, one has -- see
\cite{AK98,And02}  -- the stability result
$\|v\|_\infty\le C_i \|\phi\|_{W^{-1,\infty}}$ for a certain constant $C_i$
that depends only on $b_{ii}$ and $a_{ii}$. Apply this result to \eqref{convdiffstrong} then,
similarly to the analysis of \eqref{convdiff-itheqB}, gather the $\|u_j\|_\infty$ terms to the
left-hand side. Define the $M \times M$ matrix $ \Upsilon = (\gamma_{ij})$ by $\gamma_{ii}=1,\
\gamma_{ij} = -C_i[\|b_{ij}'+a_{ij}\|_{L_1} + \|b_{ij}\|_\infty]$ for $i\ne j$. Assuming that
$\Upsilon$ is inverse monotone,
 we get an a priori bound on $\|\vec{u}\|_\infty$.  Remark that this assumption on
 $B$ implies that $B$ is strictly diagonal dominant. In the symmetric case, then for $b_{ii}>0$
 the positive definitness of $B$ follows.

  Using this
bound, it is shown in \cite{ORSS08} that one can decompose each component of $\vec u$ into smooth
and layer components.

 The case of eigenvalues of different sign for a symmetric matrix with a special structure
is already studied in \cite{AKK74}, see also \cite{RR11} for a system of two equations.
Based on the assumption that $A+1/2B'$ is positive semidefinite, a stability result
in the energy norm allows the proof of the existence of a solution decomposition (see \cite{RR11}, Theorem 2.5).

For strongly coupled systems and different small parameters $\varepsilon_i$ the situation is even more
complicated. For first results concerning layer structure, layer location and characterization of the
reduced problem see \cite{Me11,Ro12}.

For systems of reaction-diffusion equations
\begin{subequations}\label{reacdiff}
\begin{align}
{\cal L}\vec{u}:&= -\varepsilon^2\vec{u}'' + A\vec{u}= \vec{f}
\quad\text{on }\Omega:=(0,1),\\
\vec{u}(0) &=\vec{g}_0,\ \vec{u}(1) =\vec{g}_1,
\end{align}
\end{subequations}
often the assumption
\begin{equation}\label{posdef}
\xi^TA\xi\ge \alpha^2\xi^T\xi
\end{equation}
with a constant $\alpha>0$  is used. Then, $\vec{u}$ can be decomposed into
a smooth part $\vec{v}$ and a layer part $\vec{w}$, where the layer part satisfies
\begin{equation}
|w_i^{(j)}|\le C\varepsilon^{-j}\left(e^{-dx/\varepsilon}+e^{-d(1-x)/\varepsilon}\right)
\end{equation}
with some positive constant $d$. This fact was already used by Bakhvalov \cite{Ba69} to analyse
a finite difference method for \eqref{reacdiff} on a special mesh.

In \cite{LM09} the authors study $L_\infty$ stability (assumption \eqref{posdef} implies stability in some
energy norm) of the system
\begin{subequations}\label{reacdiff2}
\begin{align}
{\cal L}\vec{u}:&= -{\cal E}\vec{u}'' + A\vec{u}= \vec{f}
\quad\text{on }\Omega:=(0,1),\\
\vec{u}(0) &=\vec{g}_0,\ \vec{u}(1) =\vec{g}_1,
\end{align}
\end{subequations}
with ${\cal E}:={\rm diag}(\varepsilon_1^2,\cdots,\varepsilon_M^2)$ and $\varepsilon_1\le\varepsilon_2\le\cdots\le\varepsilon_M$.
Assuming $a_{kk}>0$, from
$$
 -\varepsilon_k^2 u_k''+ a_{kk}u_k
    = f_k - \sum_{m\not=k}a_{km}u_m
$$
follows as above for strongly coupled systems $L_\infty$ stability if
\begin{equation}\label{gammacond}
\Gamma^{-1}\ge 0 \quad {\rm with}\quad \gamma_{ii}=1,\quad
   \gamma_{ij}=-\|\frac{a_{ij}}{a_{ii}}\|_\infty\quad {\rm for}\,\,i\not=j.
\end{equation}
Assumption \eqref{gammacond} is satisfied if $A$ is strongly diagonal dominant. In the symmetric
case then \eqref{posdef} follows.

Now let $A$ be strongly diagonal dominant and
$$
 \sum_{k\not=i}\|\frac{a_{ik}}{a_{ii}}\|_\infty<\zeta<1.
$$
Define
$$
\kappa^2:=(1-\zeta)\min_i\min_{x\in [0,1]}a_{ii}(x) \quad {\rm and}\,\,
   {\cal B}_\varepsilon (x):=e^{-\kappa x/\varepsilon}+e^{-\kappa (1-x)/\varepsilon}.
$$
Then (see \cite{LM09}, Theorem 2.4), $\vec{u}$ can be decomposed in $\vec{v}+\vec{w}$, where the
layer part satisfies
\begin{equation}\label{dec-r-d1}
|w_i^{(k)}|\le C\sum_m\varepsilon_m^{-k}{\cal B}_{\varepsilon_m} (x)\quad {\rm for}\,\,k=0,1,2
\end{equation}
and
\begin{equation}\label{dec-r-d2}
|w_i^{(k)}|\le C\varepsilon_i^{2-k}\sum_m\varepsilon_m^{-2}{\cal B}_{\varepsilon_m} (x)\quad {\rm for}\,\,k=3,4.
\end{equation}

In the case $M=2$ and analytic data assuming \eqref{posdef}, full asymptotic expansions that are explicit
in the perturbation parameters and the expansion order are presented in \cite{MXO12}.

\section{Special meshes and the analysis of FEM on DL-meshes}
\subsection{Special meshes for problems with layers}
We now present the construction of meshes suitable for a problem posed in $[0,1]$ with a layer term
$E=\exp(-\gamma\,x/\varepsilon)$ (of course, the layer could also be located at $x=1$ or at both
endpoints of the interval).

Our aim is to achieve \emph{uniform
convergence}\index{uniformly convergent scheme}  in the discrete maximum norm; that is, the
computed solution $\{u_i^N\}_{i=0}^N$ satisfies
\begin{equation}\label{eq1:unifcgce}
\|u-u^N\|_{\infty,d} := \max_{i=0,\dots,N} |u_i-u_i^N| \le C N^{-\alpha}
\end{equation}
for some positive constants $C$ and  $\alpha$ that are independent of $\varepsilon$ and of $N$.
A~power of $N$ is a suitable measure of the error $u-u^N$ for the particular families of meshes
that are discussed, but a bound of this type is inappropriate for an
arbitrary family of meshes; see \cite{SW96}.

The aim to achieve uniform convergence in the maximum norm is demanding and leads to meshes
which sometimes do not have desirable properties. Therefore, we will take into account as well
meshes where the constant $C$ in \eqref{eq1:unifcgce} will weakly depend on $\varepsilon$. Such a
desirable property is, for example, the local quasi-uniformity of the
mesh.

The simplest meshes which we discuss first are piecewise equidistant, but we also consider in detail
graded meshes. Because for $x=\mu \frac{\varepsilon}{\gamma}\ln  \frac{1}{\varepsilon}$ the layer
term is of order $O(\varepsilon^\mu)$, a simple idea is to use a fine equidistant mesh in the
interval $[0,\mu \frac{\varepsilon}{\gamma}\ln  \frac{1}{\varepsilon}]$ and a coarse
equidistant mesh in $[\mu \frac{\varepsilon}{\gamma}\ln  \frac{1}{\varepsilon},1]$. But with such
an so called A-mesh, in general, it is not possible to achieve uniform convergence.

Shishkin spread the idea to choose  as a transition point $\sigma$ from the fine to the coarse mesh
the point defined by
\begin{equation}\label{transition point S}
\sigma=\min\{1/2,\mu \frac{\varepsilon}{\gamma}\ln  N \}
\end{equation}
and to subdivide each of the intervals $[0,\sigma]$ and $[\sigma,1]$ by an equidistant mesh with
$N/2$ subintervals. Then, for small $\varepsilon$, one has $E(\sigma)=N^{-\mu}$, and $\mu$ is
chosen in dependence of the order of the method used.

\begin{remark}\label{rem:shishpointdist}
It is not vital  that one has exactly the same number of
subintervals in $[0,\sigma]$ and $[\sigma,1]$. All that the
theory demands is that as $N\rightarrow\infty$ the number of
subintervals in each of these two intervals is bounded below by
$CN$ for some constant $C>0$. \hfill$\clubsuit$
\end{remark}

The coarse  part of this Shishkin mesh has spacing
$H=2(1-\sigma)/N$, so $N^{-1}\le H \le 2N^{-1}$. The fine part has
spacing $h=2\sigma/N = (4/\gamma)\varepsilon N^{-1}\ln N$, so
$h\ll\varepsilon$. Thus there is a very abrupt change in mesh size
as one passes from the coarse part to the fine part. The mesh is
not locally quasi-equidistant,  uniformly in $\varepsilon$.

\begin{remark}\label{rem:ShishDecay}
\emph{(A key property of the Shishkin mesh)} Nonequidistant  meshes for convection-diffusion
problems are sometimes described as \lq\lq layer-resolving" meshes. One might infer from this
terminology that wherever the derivatives of $u$ are large, the mesh is chosen so fine that the
truncation error of the difference scheme is controlled. But the Shishkin mesh does not fully resolve the layer: for
$$|u'(x)| \approx C\varepsilon^{-1}\exp(-\gamma\,x/\varepsilon),$$
so
$$|u'(\sigma)| \approx C\varepsilon^{-1}\exp(-\mu\ln N) = C\varepsilon^{-1}N^{-\mu},$$
which can be large.
 That is, $|u'(x)|$ can still be large on part of the first
coarse-mesh interval $[x_{N/2}, x_{N/2+1}]$.

At first sight this incomplete resolution of the boundary layer seems like a flaw, but it is in
fact the key property of the mesh! Shishkin's insight was that one could achieve satisfactory
theoretical and numerical results without resolving all of the layer and as a consequence his mesh
permits us to use a fixed number of mesh points that is independent of $\varepsilon$. If one set
out to \lq\lq repair" the Shishkin mesh  by constructing a two-stage
piecewise-equidistant mesh as we have done, but with the additional requirement that the mesh be
fine enough to control the local truncation  error wherever the derivatives of $u$ are very large,
then the number of mesh points required would have to grow like $|\ln\varepsilon|$ as
$\varepsilon\rightarrow 0$.  See \cite{Wes96} and \cite[Section 3.6]{FHMOS00}.

Although the number of mesh points is fixed independently of $\varepsilon$, nevertheless numerical
analysis on Shishkin meshes does pay a price for the nature of the
construction: typically the trickiest part of the domain to handle is the first coarse mesh
interval -- because the derivatives of $u$ can be large there. \hfill$\clubsuit$
\end{remark}

If a method for a problem with a smooth solution has the order $\alpha$, due to the fine mesh
size $h=O(\varepsilon N^{-1}\ln N)$ in the case $u^{(k)}\approx \varepsilon ^{-k}$ we can expect
 that the error on a Shishkin mesh is of the order $O((N^{-1}\ln N)^\alpha)$. Especially for higher
order methods the logarithmic factor is not nice. An optimal mesh should generate an error of the
order $O(N^{-\alpha})$.

One possibility to generate an optimal mesh is to introduce a graded mesh in the fine subinterval
$[0,\sigma]$ using a {\it mesh-generating function}. Assuming the
function
$\lambda:\,\, [0,1/2]\mapsto [0, \ln N]$ to be strictly increasing, set
$$
   x_i= \frac{\mu\varepsilon}{\gamma}\lambda(i/N),\quad i=0,1,\cdots, N/2.
$$
We call such meshes {\it Shishkin-type meshes}\index{Shishkin-type mesh}. It turns out (see the
next Section) that in error estimates for Shishkin-type meshes often the factor $\max|\psi'(\cdot)|$
appears, here $\psi$ is the {\it mesh-characterizing function} defined by
$$
   \psi:=e^{-\lambda}:\quad [0,1/2]\mapsto [1, 1/N].
$$
For the original Shishkin mesh we have $\max|\psi'(\cdot)|=O(\ln N)$. A popular optimal mesh
is the {\it Bakhvalov-Shishkin mesh}\index{Bakhvalov-Shishkin mesh} with
$$
 \psi(t)=1-2t(1-N^{-1}) \quad {\rm and}\quad \max|\psi'(\cdot)|\le 2.
$$
The mesh points of the fine mesh are given by
\begin{equation}\label{BS-mesh}
 x_i= -\frac{\mu\varepsilon}{\gamma}\ln\left(1-2(1-N^{-1})\frac{i}{N}\right),\quad i=0,1,\cdots, N/2.
\end{equation}
For other possibilities to choose $\lambda$, see \cite{RL99}.

\begin{remark}
In \cite{FX17} we find a slight generalization of Shishkin-type meshes based on the property
$\lambda(1/2)=\ln (\theta N)$ with some additional parameter $\theta$. This allows to
characterize the so called eXp-mesh as generalized Shishkin-type mesh.
\hfill$\clubsuit$
\end{remark}

If one chooses the transition point from a fine to the coarse mesh by
\begin{equation}\label{transition point B}
\sigma^*=\min\{1/2,\mu \frac{\varepsilon}{\gamma}\ln \frac{1}{\varepsilon}\},
\end{equation}
one has to use a graded fine mesh to achieve uniform convergence. Similarly as \eqref{BS-mesh}, a
{\it Bakhvalov-type mesh}\index{ Bakhvalov-type mesh} is given by
\begin{equation}\label{B-type-mesh}
 x_i= -\frac{\mu\varepsilon}{\gamma}\ln\left(1-2(1-\varepsilon)\frac{i}{N}\right),\quad i=0,1,\cdots, N/2.
\end{equation}
In $[\sigma^*,1]$ the mesh is equidistant.

Bakhvalov's original mesh \cite{Ba69}\index{Bakhvalov mesh} was a little more complicated. The mesh points
near $x=0$ were defined by
$$
    q(1-\exp(-\frac{\gamma x_i}{\mu\varepsilon}))=\frac{i}{N}
$$
 in some interval $[0,\tau]$, here $q$ is a parameter. That means, Bakhvalov proposed the mesh-generating function
$$
\phi(t)=\left\{ \begin{array}{ll}
    -\frac{\mu\varepsilon}{\gamma}\ln \frac{q-t}{t}\quad {\rm for}\quad t\in[0,\tau]\\
    \phi(\tau)+\phi'(\tau)(t-\tau)\quad {\rm for}\quad t\in[\tau,1].
    \end{array}\right.
$$
Remarkably,  $\tau$ is defined by the requirement that the mesh-generating function is $C^1$.
Thus, $\tau$ has to solve the nonlinear equation
$$
   \phi'(\tau)=\frac{1-\phi(\tau)}{1-\tau}.
$$
Bakhvalov-type meshes are simpler than the Bakhvalov meshes, the mesh-generating function is not
longer $C^1$. But both meshes are not locally quasi-equidistant. In some cases for these meshes
optimal error estimates are known, but the analysis is often more complicated than for
Shishkin-type meshes. In some papers Vulanovi\'c \cite{Vu83,Vu86,Vu89,Vu91} simplified the
mesh-generating function of a Bakhvalov mesh such that the resulting equation for $\tau$ is
easy to solve.

As Lin{\ss} pointed out \cite{Li10}, a Bakhvalov-mesh can also be generated by equidistributing
the monitor function
$$
M(s)=\max\left(1, \tilde K\gamma\varepsilon^{-1}e^{-\frac{\gamma s}{\mu\varepsilon}}\right).
$$

\begin{remark}
In several papers (see \cite{LP89} and its references)  Liseikin examines the convergence of
finite difference methods when using mesh generating functions $\lambda(t)$ of the given
independent variable that satisfy $|\lambda'(t)| \le C$ for all $t\in [0,1]$. This approach
generates a graded grid of Bakhvalov type\index{Bakhvalov-type mesh}. His book \cite{Lis99}
develops a general theory of grid generation. The analysis in these sources is written in terms of
``layer-resolving transformations"; their relationship to mesh generating functions in a singular
perturbation context is discussed in \cite{Vul07}.
\hfill$\clubsuit$
\end{remark}
\begin{remark}
In \cite{SY12} we find the proposal to generate a mesh by the implicitly defined function
\begin{equation}\label{Lambert}
 \xi(t)-e^{\frac{\gamma\xi(t)}{\mu\varepsilon}}+1-2t=0.
\end{equation}
The mesh has the advantage that it is not necessary to use different mesh generating functions in
different regions. But \eqref{Lambert} is not so easy to solve, however, a solution based on the
use of Lambert's W-function is possible. Some difference schemes (and finite elements) on that
mesh can be analyzed similarly as on Bakhvalov-type meshes \cite{RTU14}.
\hfill$\clubsuit$
\end{remark}

So far we constructed a graded fine mesh using a mesh-generating function. Alternatively, it
is also relatively popular to use a {\it recursive formula}.

Gartland \cite{Ga88} graded a mesh in the following way:
$$
x_0=0,\quad x_1=\varepsilon H,\quad x_{i+1}=x_i+h_i
$$
with
\begin{equation}\label{Gartland}
  h_i=\min\left(H, \varepsilon H e^{\frac{\gamma x_i}{2\varepsilon}},e h_{i-1}.
   \right)
\end{equation}
The restriction $h_i\le e h_{i-1}$ ensures that the mesh is locally quasi-equidistant.

\begin{remark}\label{rem1:upwindarbmesh}
If simple upwinding \index{simple upwind scheme} for a convection-diffusion problem is uniformly
convergent \index{uniformly convergent scheme} in the sense of \eqref{eq1:unifcgce} for some
constant $\alpha >0$, and the mesh is locally quasi-equidistant \index{locally quasi-equidistant
grid} (uniformly in $\varepsilon$), then the number $N$ of mesh intervals must increase as
$\varepsilon\rightarrow 0$. To see this, observe that the arguments of \cite{St03} are still valid
when slightly modified by considering a limit as $N\rightarrow\infty$ with $\varepsilon \ge h_1$
and $i=1$; one then arrives at the conclusion of that paper that $h_1 = o(\varepsilon)$. (There
are some minor extra mesh assumptions such as existence of $\lim_{N\rightarrow\infty} h_1/h_2$ and
$\lim_{N\rightarrow\infty} h_2/h_1$.) But the mesh diameter is at least $1/N$, so the locally
quasi-equidistant property implies that $\varepsilon K^N \ge 1/N$, where $K$ is the constant in
$$
   h_i\le Kh_j \quad {\rm for}\quad |i-j|\le 1.
$$
 Hence $NK^N \ge 1/\varepsilon$, so $N \approx \log_K(1/\varepsilon)$.
 \hfill$\clubsuit$
\end{remark}

Introducing the transition points $x^*,x'$ by
$$
    x^*\approx K\varepsilon \ln \frac{K}{H}, \quad x'\approx K\varepsilon \ln \frac{K}{\varepsilon}
$$
Gartland observed that the number of mesh points in the inner region $[0,x^*]$ as well in the
outer region $[x',1]$ is of order $O(1/H)$, but in the transition region $[x^*,x']$ of order
$O(\ln\ln \frac{H}{\varepsilon})$.

We call the modification of the mesh where \eqref{Gartland} is replaced by
\begin{equation}
 h_i=\min\left(H, \varepsilon H e^{\frac{\gamma x_i}{2\varepsilon}}.
   \right)
   \end{equation}
{\it Gartland-type mesh}\index{Gartland-type mesh}. The number of mesh points is now independent
of $\varepsilon$ and the mesh allows optimal error estimates. The mesh is not locally
quasi-eqidistant.
 Gartland-type meshes are studied
in \cite{RS97}, \cite{LX06,LX09,LWX16}.

Much simpler is the {\it Duran-Lombardi mesh} \cite{DL06}\index{Duran-Lombardi mesh} defined by
$$
\begin{array}{ll}
x_0=0,\,\,x_i=i\kappa H\varepsilon\quad {\rm for}\quad 1\le i\le \frac{1}{\kappa H}+1\\
 x_{i+1}=x_i+\kappa H x_i\varepsilon\quad {\rm for}\quad  \frac{1}{\kappa H}+1\le i\le M-2,\quad x_M=1.
\end{array}
$$
Here M is chosen such that $x_{M-1}<1$ but $x_{M-1}+\kappa H x_{M-1}\ge 1$, assuming that the last interval
is not extremely small.

The mesh is locally quasi-equidistant and glitters by its simplicity, there is no need to define a transition
point. Uniform error estimates with respect to $H$ are possible, but the number of mesh-points is
proportional to $\frac{1}{H}\ln \frac{1}{\varepsilon}$. Numerical experiments show that the mesh is more
robust than other meshes with respect to the use of a mesh designed for some given $\varepsilon$ for
larger values of the parameter.

If the mesh is given by a recursive formula, one can define a mesh-generating function by interpolation
of the values given in the mesh points. That gives also a possibility to analyze discretization
methods on such a mesh, see, for instance, \cite{RTU15}.

\begin{remark}
When analyzing hp finite element methods for singularly perturbed problems it is common to use
an hp boundary layer mesh, see \cite{Me02}. For such methods it is possible to prove
exponential convergence.
\hfill$\clubsuit$
\end{remark}

\begin{remark}
For problems with weak layers it is under some conditions possible to use a coarser mesh in comparison to meshes used for strong layers, see \cite{Ro22}.
\hfill$\clubsuit$
\end{remark}

\subsection{The analysis of finite element methods on DL meshes}

The analysis of finite element methods for standard convection-diffusion problems with
exponential boundary layers for Shishkin-type meshes is well known. One obtains for
linear or bilinear elements
\begin{equation}
\|u-u^N\|_\varepsilon\le CN^{-1}\max|\psi'|
\end{equation}
with the mesh characterizing function $\psi$ introduced in the previous subsection. Consequently, a Bakhvalov-Shishkin
mesh or a Gartland-Shishkin mesh yield the optimal estimate
\begin{equation}
\|u-u^N\|_\varepsilon\le CN^{-1}.
\end{equation}

For Bakhvalov meshes the proof of an optimal error estimate in 2D was for a long time an open problem, the analysis in
Roos/Schopf \cite{RS12}  yield only a weak dependence on $\varepsilon$. But beginning in 2020 Jin Zhang and Xiaowei Liu published
a new approach yielding optimal estimates using a new tricky interpolant, see \cite{ZL20}, \cite{ZL21}.

Next we sketch as an example for the analysis on recursively defined meshes the analysis on Duran-Lombardi meshes.
For the analysis presented it is sufficient to assume \\
for $0\le k\le 2$
\begin{align*}
|\frac{\partial^k u}{\partial x^k}|&\le C(1+\varepsilon^{-k}\exp(-\beta_1 x/\varepsilon)), \\
|\frac{\partial^k u}{\partial y^k}|&\le C(1+\varepsilon^{-k}\exp(-\beta_2 y/\varepsilon))
\end{align*}
and
\begin{equation*}
|\frac{\partial^2 u}{\partial x\partial y}|\le C(1+\varepsilon^{-1}\exp(-\beta_1 x/\varepsilon)+
  \varepsilon^{-1}\exp(-\beta_2 y/\varepsilon)+\varepsilon^{-2}\exp(-\beta_1 x/\varepsilon)\exp(-\beta_2 y/\varepsilon)).
\end{equation*}
(for that analysis we consider layers located at $x=0$ and at $y=0$)

Consider bilinear elements on a mesh defined by
$$
\begin{array}{ll}
x_0=0,\,\,x_1=\kappa H\varepsilon\\
 x_{i+1}=x_i+\kappa H x_i\quad {\rm for}\quad  1\le i\le M-2,\quad x_M=1.
\end{array}
$$
Here M is chosen such that $x_{M-1}<1$ but $x_{M-1}+\kappa H x_{M-1}\ge 1$, assuming that the last interval
is not extremely small. Remark that this mesh was the first proposal of Duran and Lombardi, the modification
presented is characterized by an increasing mesh size and can be analyzed similarly. As mentioned before
the number of mesh points $N$ is proportional to $H^{-1}|\ln \varepsilon|$. Analogously we define
the mesh points $y_j$ in the $y$ direction.

First one obtains for the interpolation error on that mesh
\begin{equation*}
\|u-u^I\|_0\le CH^2 \quad{\rm and}\,\,\varepsilon^{1/2}|u-u^I|_1\le CH.
\end{equation*}
Let us sketch, for instance, the proof of the estimate in the $H^1$ semi norm. Consider an element
$K_{ij}=(x_{i-1},x_i)\times(y_{j-1},y_j)$. Then the anisotropic interpolation error estimate reads
$$
\|(u-u^I)_x\|_{0,K_{ij}}^2\le C\left(h_i^2\|u_{xx}\|_{0,K_{ij}}^2+
    h_j^2\|u_{yy}\|_{0,K_{ij}}^2\right).
$$
On $K_{11}$ the desired estimate follows immediately because $h_1=k_1=H\varepsilon$.\\
Consider next the case $i,j\ge 2$. Then, the use of the recursive formula $h_i=\kappa Hx_{i-1}$
allows the estimate
$$
\|(u-u^I)_x\|_{0,K_{ij}}^2\le CH^2\left(\|xu_{xx}\|_{0,K_{ij}}^2+
    \|yu_{yy}\|_{0,K_{ij}}^2\right).
$$
Now the analysis is based on the simple observation that
$$
\int_0^1e^{-x/\varepsilon}=O(\varepsilon)\quad{\rm but}\,\,\int_0^1xe^{-x/\varepsilon}=O(\varepsilon^2)
 \quad{\rm and}\,\,\int_0^1x^{1/2}e^{-x/\varepsilon}=O(\varepsilon^{3/2}).
$$
Using that observation our estimates for the derivatives of $u$ yield, for instance,
$$
\|u_{xx}\|_0\le C\varepsilon^{-3/2} \quad{\rm and}\quad\|xu_{xx}\|_0\le C\varepsilon^{-1/2},
$$
 similarly
$$
\|u_{xy}\|_0\le C\varepsilon^{-1} \quad{\rm and}\quad\|y^{1/2}u_{xy}\|_0\le C\varepsilon^{-1/2}.
$$
From these estimates the desired estimate for $\|(u-u^I)_x\|_0$ follows.\\
The remaining elements $K_{ij}$ can be handled similarly.

The error analysis follows the same lines as on Shishkin meshes, the only critical part is to estimate
$$
  I_{conv}=\int_{\Omega}b\nabla(u-u^I)(u^I-u_H).
$$
Introducing the subdomain $\Omega_3$ by
$$
   \bar \Omega_3=\cup\left\{ \bar K_{ij}:\quad x_{i-1}\ge c\varepsilon|\ln \varepsilon|,\quad
    y_{j-1}\ge c\varepsilon|\ln \varepsilon|\right\}
$$
we get on $\Omega_3$
$$
  |\int_{\Omega_3}b\nabla(u-u^I)(u^I-u_H)|\le CH\|u^I-u_h\|_0
$$
because the behavior of $u$ implies $|u-u^I|_{1,\Omega_3}\le CH$.

Consider next, for example, the strip
$$
   \bar \Omega_1=\cup\left\{ \bar K_{ij}:\quad x_{i-1}< c\varepsilon|\ln \varepsilon|
    \right\}.
$$
The Poincar\'e-Friedrichs inequality
$$
  \|u^I-u_H\|_{0,\Omega_1}\le C\varepsilon \ln\frac{1}{\varepsilon}\|\frac{\partial}{\partial x}(u^I-u_H)\|_{0,\Omega_1}
$$
yields
$$
  |\int_{\Omega_1}b\nabla(u-u^I)(u^I-u_H)|\le C\varepsilon^{1/2}\|u^I-u\|_1 \ln\frac{1}{\varepsilon}
       \|u^I-u_H\|_\varepsilon,
$$
resulting finally in the error estimate
$$
   \|u-u_H\|_\varepsilon \le CH\ln\frac{1}{\varepsilon}.
$$
Introducing the number of mesh points used, one has
$$
   \|u-u_H\|_\varepsilon \le CN^{-1}(\ln\frac{1}{\varepsilon})^2.
$$
This estimate on the Duran-Lombardi mesh is not uniform in $\varepsilon.$ But the analysis requires
only estimates for derivatives, and two further advantages of the mesh are its simplicity (no definition
of a transition point is necessary) and its robustness (a mesh defined for some $\varepsilon$ can also be used
for larger values of the parameter).

2021 Y. Yin and P. Zhu found a further remarkable property of D-L-meshes: for the streamline-diffusion finite element method
one can prove uniform estimates also in the full streamline-diffusion norm, see \cite{YZ19} and \cite{BR21}.

\section{Numerical Methods for Systems}\label{sec1:2.6}
Consider again the system
\begin{subequations}\label{convreacdiff2}
\begin{align}
L\vec{u}:&= -\varepsilon\vec{u}''+B\vec{u}' + A\vec{u}= \vec{f}
\quad\text{on }\Omega:=(0,1),\\
\vec{u}(0) &=\vec{g}_0,\ \vec{u}(1) =\vec{g}_1.
\end{align}
\end{subequations}
Based on the analytical results presented in Section 2 we now study numerical methods, extending some
results from the survey \cite{LS09}.

In \cite{AKK74} the authors consider simple and midpoint upwind on a uniform mesh in a case of strong
coupling, assuming $B$ is symmetric and has the special structure
$$
B=\left(
\begin{array}{ll}  B^1\quad 0\\
                  0 \quad  B^2
\end{array}
\right).
$$
Here, for instance, $B^1$ is positive definite, $B^2$ negative definite. Far from the layers error estimates
are proven, but,
of course, one cannot obtain uniform convergence.

In \cite{HYY17} the Il'in-Allen-Southwell scheme is generalized for systems. Let us start from
$$
 -\varepsilon\vec{u}''(x_i)+B(x_i)\vec{u}'(x_i) + A(x_i)\vec{u}(x_i)= \vec{f}(x_i)
$$
and assume the matrix $B(x_i)$ to be symmetric with eigenvalues $\lambda_j^i$. Then, there exists a matrix
$P_i$ which diagonalizes $B(x_i)$. Setting
$$
     \vec{v}=P_i^{-1} \vec{u},\quad \vec{g}=P_i^{-1} \vec{f}
$$
we obtain for the j-th component of $\vec{v}$ the equation
$$
 -\varepsilon v_j''(x_i)+\lambda_j^i v_j'(x_i) + \sum_k a^*_{jk}v_k(x_i)= g_j(x_i).
$$
Next we introduce the scalar version of the  Il'in-Allen-Southwell scheme:
$$
 -\varepsilon \sigma_j^i D^+D^-v_j^i+\lambda_j^i D^0 v_j^i + \sum_k a^*_{jk}v_k^i= g_j(x_i)
$$
with $\sigma_j^i=\frac{\lambda_j^i\,h}{2}\coth \frac{\lambda_j^i\,h}{2\varepsilon}$. Collecting the
$\sigma_j^i$ in a vector, the back transformation yields
\begin{equation}\label{IAS}
 -\varepsilon P_i\sigma^iP_i^{-1}\vec{u^i}+B(x_i)D^0\vec{u^i} + A(x_i)\vec{u^i}= \vec{f}(x_i).
\end{equation}
If $B$ is not diagonalizable one can use the Jordan canonical form, the resulting scheme is a little
more complicated. The numerical experiments in \cite{HYY17} promise first order uniform convergence,
but there exists no analysis.

Next we sketch error estimates for layer-adapted meshes, mostly for Shishkin meshes. Of course, to create these meshes we need
a priori information on the location of layers and precise estimates on the behavior of derivatives,
which we discussed for systems in Section 2.

Under the assumptions formulated in that Section for weakly coupled systems, even with different diffusion
coefficients, one can analyze the upwind finite difference scheme based on the information for the first
order derivatives given in \eqref{syst,deriv}. By the transformation $x:=1-x$ we shift the layer to
$x=0$. Then, the mesh is equidistant and coarse away from $x=0$, and piecewise equidistant with successively
finer meshes, as one approaches $x=0$. More precisely, we have mesh transition points $\tau_k$ defined by
\begin{equation}\label{sys,mesh}
\tau_{M+1}=1,\quad \tau_k=\min \left(\frac{k\tau_{k+1}}{k+1},\frac{\sigma\varepsilon_{k}}{\beta}\ln N
\right)\quad {\rm and}\quad
   \tau_{0}=0.
\end{equation}
Then, for $k=0,\cdots,M$, the mesh is obtained by dividing each of the intervals $[\tau_k,\tau_{k+1}]$
into $N/(M+1)$ subintervals of equal length.

Similarly as in the scalar case now from a $(\|\cdot\|_{\infty,d},\|\cdot\|_{-1,\infty,d})$ stability result
one can prove for simple upwind
$$
  \|\vec{u}-\vec{u^N}\|_{\infty,d}\le C\max_k \int_{x_{k-1}}^{x_k}(1+|\vec{u'}|)dx,
$$
and
$$
  \|\vec{u}-\vec{u^N}\|_{\infty,d}\le C N^{-1}\ln N \quad {\rm for\,\, Shishkin\,\, meshes}
$$
and first order uniform convergence for Bakhvalov meshes follows \cite{Li07}.

Results in the maximum norm for higher order schemes are not known. In \cite{RS15} linear finite elements
are analyzed in the energy norm on a modified Shishkin mesh, but also studied numerically in the maximum
norm. It seems that the weak layers for higher order derivatives reduce the convergence order two in
the maximum norm observed on the modified mesh if the standard Shishkin mesh is used instead.

For strongly coupled system we assume $\varepsilon_i=\varepsilon$ for all $i$. In the case of
overlapping layers at one boundary with an analogous technique as for weakly coupled systems
one can prove for the upwind scheme on a Shishkin mesh
$$
  \|\vec{u}-\vec{u^N}\|_{\infty,d}\le C N^{-1}\ln N,
$$
see \cite{ORSS08}. In the case $M=2$, where both solution components do have layers at $x=0$ and
at $x=1$, linear finite elements are analyzed in \cite{RR11} in some energy norm.

Consider finally a system of reaction-diffusion equations:
\begin{subequations}
\begin{align}
{\cal L}\vec{u}:&= -{\cal E}\vec{u}'' + A\vec{u}= \vec{f}
\quad\text{on }\Omega:=(0,1),\\
\vec{u}(0) &=\vec{g}_0,\ \vec{u}(1) =\vec{g}_1,
\end{align}
\end{subequations}
On a uniform mesh, in \cite{Li16} a method is proposed which produces solutions
without numerical oscillations. The method starts with a least-squares functional
$$
 {\cal F}(\vec{u},\vec{p})=\|D((\vec{u},\vec{p}))-F\|_0^2
$$
where $F=[0,f]^T$ and
$$
D((\vec{u},\vec{p}))=[\varepsilon_1(p_1-\nabla u_1),\cdots,\varepsilon_M(p_M-\nabla u_M),
    -{\cal E}\nabla\vec{p} + A\vec{u}]^T.
$$
A dG method is used to discretize the least-square problem with linear finite elements, see
\cite{Li16}.

Under the assumptions formulated in Section 2 we have bounds for derivatives of the solution
and can construct a layer adapted mesh. Consider the Shishkin mesh of the type \eqref{sys,mesh},
now we use that mesh at $x=0$ and at $x=1$.

It is well known that  $L_\infty$ stability is not sufficient to analyze the standard (central)
finite difference scheme on a Shishkin mesh. For a system of two equations in \cite{LM04} the
maximum principle and special barrier functions are used to prove
$$
  \|\vec{u}-\vec{u^N}\|_{\infty,d}\le C (N^{-1}\ln N)^2.
$$
Later in the general case \cite{LM09} the authors started again with the discrete $L_\infty$ stability, but used next
a representation of the consistency error via the discrete Green's function:
$$
\eta_{i,j}=\sum_k \bar h_kG_k\varepsilon_i^2(D^+D^-u_i-u''_i).
$$
A tricky manipulation of this representation which uses the properties of the discrete Green's function and the
bounds \eqref{dec-r-d1}, \eqref{dec-r-d2} leads to
$$
\|\vec{u}-\vec{u^N}\|_{\infty,d}\le C\left(
  \max_k \int_{x_{k-1}}^{x_k}(1+\sum_m \varepsilon_m^{-1}B_{2\varepsilon_{m}}(x))dx
\right)^2.
$$
On a Shishkin mesh follows
$$
  \|\vec{u}-\vec{u^N}\|_{\infty,d}\le C (N^{-1}\ln N)^2.
$$
It is also possible to obtain the corresponding slightly better results on Bakhvalov meshes.

Linear finite elements on layer-adapted meshes in the energy norm are analyzed in \cite{Li09}, robust
exponential convergence of hp FEM for a set of two equations in \cite{MXO13}. So far error estimates
for systems in a balanced norm are rare.

\section{Remarks to the discontinuous Galerkin method}
In \cite{RST08} we presented a short introduction into the discontinuous Galerkin method because during that
time the method was still not thus popular as today. Since then
 five books on dGFEM were published which proves the attractivity of the method
\cite{HW08,Ka08,Ri08,DE12,DF15}.

We discussed 2008 mainly the primal formulation of dGFEM and, especially, the symmetric and non-symmetric
method with interior penalties, SIP and NIP. In this Section we shortly introduce the very popular
local discontinuous Galerkin method LDG, hybrid methods HdG and comment discontinuous Petrov-Galerkin
methods dPG with optimal test functions.

\subsection{The local discontinuous Galerkin method LDG}
Consider the convection-diffusion problem
\begin{subequations}\label{cd}
  \begin{align}
Lu:=-\varepsilon \Delta u -b\cdot\nabla u+cu&=f
\qquad \mbox{in}\,\, \Omega\\
u&=0 \qquad \mbox{on}\,\, \Gamma= \partial \Omega.
   \end{align}
\end{subequations}
 Alternatively to the primal formulation one has the flux formulation of the dGFEM, which
starts from the formulation of \eqref{cd} as
\[
\theta=\nabla u,\qquad -\varepsilon\nabla\cdot\theta+cu=f.
\]
A corresponding weak form is
\begin{align*}
\int_\kappa \theta\cdot\tau &=-\int_\kappa
        u\,\,\nabla\cdot\tau+\int_{\partial\kappa}
        u\,\mu_{\kappa}\cdot\tau,\\
\varepsilon\int_\kappa \theta\cdot\nabla v+\int_\kappa
cuv&=\int_\kappa
        fv+\int_{\partial\kappa}
        \theta\cdot\mu_{\kappa}\,v.
\end{align*}
This generates the following discretization: find $u_h,\,\theta_h$ such that
\begin{align*}
\int_\kappa \theta_h\cdot\tau_h &=-\int_\kappa
        u_h\,\,\nabla\cdot\tau_h+\int_{\partial\kappa}
        \hat u_\kappa\,\mu_{\kappa}\cdot\tau_h,\\
\varepsilon\int_\kappa \theta_h\cdot\nabla v_h+\int_\kappa
cu_hv_h&=\int_\kappa
        fv_h+\int_{\partial\kappa}
        \hat\theta_\kappa\cdot\mu_{\kappa}\,v_h.
\end{align*}
Here the choice of the numerical fluxes $\hat\theta_\kappa$ and $\hat u_\kappa$ that approximate
$\theta=\nabla u$ and $u$ on $\partial\kappa $ is very important. In \cite{ABC02} one finds a
thorough discussion of 9 variants of the dGFEM that are characterized by different choices of
$\hat\theta_\kappa$ and $\hat u_\kappa$. For each of these methods, the properties of the
associated  primal formulation that is obtained by eliminating $\theta_h$ are discussed.

Now the very popular local discontinuous Galerkin method is characterized by the following choice of the numerical fluxes on interior edges:
$$
\begin{array}{ll}
 \hat u_\kappa=     \left\{ u_h \right\}+C_{12}[u_h],\\
 \hat\theta_\kappa= \left\{\theta_h\right\}-C_{12}[\theta_h]-C_{11}[u_h].
\end{array}
$$
Here $C_{12}\ge 0$, and an adequate choice of $C_{11}$ ensures stability of the method.
See \cite{DE12}, Chapter 4.4 or \cite{DF15}, Chapter 3.3 for details and \cite{ZZ12,ZZ14} for
the analysis of LDG for convection-diffusion problems on layer-adapted meshes. For time-dependent problems, see \cite{CM22}.

\subsection{Hybrid dG (HdG)}

Introducing again $\theta=\nabla u$, we start from the mixed formulation
\begin{align*}
\int_\Omega \theta\cdot\tau +\int_\Omega u {\rm div}( \tau)&=0\quad {\rm for\,\,all}\,\,\tau\in H(div,\Omega),\\
-\varepsilon\int_\Omega {\rm div}(\theta) v+\int_\Omega cuv&=(f,v)\quad {\rm for\,\,all}\,\,v\in L^2(\Omega).
\end{align*}
Now, instead of requiring the discrete fluxes to be in $ H(div,\Omega)$, completely discontinuous elements can
be used ensuring the continuity of the normal fluxes over element interfaces by adding appropriate
constraints. The corresponding discrete problem reads:\\
Find $(\theta_h, u_h, \lambda_h) $ ($\lambda_h$ lives on the set of faces), such that
\begin{align*}
\int_\Omega \theta_h\cdot\tau_h +\int_\Omega u_h {\rm div}( \tau_h)+
       \langle \lambda_h,\tau_h\cdot n\rangle_{\partial {\cal T}_h}&=0,\\
-\varepsilon\int_\Omega {\rm div}(\theta_h) v_h+\int_\Omega cu_hv_h&=(f,v),\\
        \langle \theta_h\cdot n,\mu_h\rangle_{\partial {\cal T}_h}&=0.
\end{align*}
This method, called MH-dG \cite{ES10}, can alternatively be generated by the HdG method of \cite{FQZ15}
choosing the stabilization parameter used in that method suitably, see also \cite{CGL09,Ji16}. In
\cite{CSB05} the name discontinuous Petrov-Galerkin method was used for a hybrid technique, today
this name is mostly related to a hybrid method with optimal test functions (see next
subsection ).

Let us finally notice that the {\it weak Galerkin method} (WG)\index{weak Galerkin method (WG)} is
closely related to the HdG framework \cite{CFX17}.

\subsection{Discontinuous Petrov-Galerkin methods (dPG) with optimal test functions}

Demkowicz and Gopalakrishnan developed in the years 2009-2011 a general finite element frame work
which allows to combine discontinuous Galerkin methods with the early methodology proposed by
Morton (mentioned in Section 2.2.2 of the book in 2008) of optimal test functions. First pure convection was studied
\cite{DG10}, then beginning with \cite{CHBD14} and several other papers \cite{NCC13,CEQ14} also
convection-diffusion was taken into consideration. For reaction-diffusion in the singularly perturbed case, see
\cite{HK17}.

Next we sketch the basic philosophy following the survey \cite{DG15}. Assume that we want to solve the
following problem:\\
Find $u\in U$ such that
\begin{equation}\label{hdG1}
 A(u,v)=f(v)\quad {\rm for \,\,all}\quad v\in V.
\end{equation}
We propose that the assumptions of the Babuska-Brezzi theory are satisfied. If $U_h\subset U$ is  some
finite element space, we define the test space of a Petrov-Galerkin method by $\Theta(U_h)$, where the
trial-to-test operator $\Theta: U\mapsto V$ is given by
\begin{equation}\label{hdG2}
 \langle \Theta  u,v\rangle=A(u,v)\quad {\rm for \,\,all}\quad v\in V.
\end{equation}
Here $\langle \cdot \rangle$ denotes a scalar product in $V$. Then, the Petrov-Galerkin method
\begin{equation}\label{hdG3}
 A(u_h,v_h)=f(v_h)\quad {\rm for \,\,all}\quad v_h\in \Theta(U_h)
\end{equation}
with $u_h\in U_h$ has the optimal property
$$
\|u-u_h\|_E=\inf_{w_h\in U_h}\|u-w_h\|_E
$$
with respect to the norm
$$
   \|w\|_E:=\sup \frac{A(w,v)}{\|v\|_V}.
$$
The success of this strategy depends on how easily one can compute the test functions. To get
a local, element-by-element computation discontinuous methods are used. But still there exist many
possibilities of, in general, hybrid formulations and choices of the broken Sobolev spaces spaces $U$ and $V$ and
corresponding broken norms \cite{DG15}. So far we described the {\it ideal} dPG method, practically one
has to approximate the test functions. In the singularly perturbed case it is open to realize that
in a robust way.

Remark that there also exist close relations of ideal dPG to a method of Cohen, Dahmen, Huang, Schwab and Welper
\cite{CDW12,DHSW12}.

\section{Adaptive methods}
\subsection{The stationary case}
In the book 2008 we described briefly four types of estimators:
\begin{itemize}
\item residual estimators
 \item  estimators based on the solution of local problems
 \item  estimators based on averaging/pre-processing of the flux
 \item goal-oriented estimators (or the DWR method\index{DWR method}: dual weighted residuals).
 \index{dual weighted residuals (DWR) method}
\end{itemize}

In 2008,  estimators based on averaging/pre-processing of the flux were still not thus popular,
let us shortly introduce these estimators.
 The starting
point is the introduction of some $q\in H({\rm div};\Omega)$ in the residual equation:
\begin{align}
 {\rm Res}(v)=(f+{\rm div}\,q,v)-(\nabla u_h-q,\nabla v)\quad
        \forall v\in V. \label{flux}
\end{align}
In general, $q$ is designed in such a way that it satisfies the {\it equilibration condition}
\begin{equation}
{\rm div}\,q+Pf=0
\end{equation}
where $Pf$ is some projection of $f$. Then, for non-singularly perturbed problems the first term
in \eqref{flux} is small and the second term generates a good estimator depending on the
concrete realization of $q$. In \cite{CM10} the authors present five methods to define $q$ and
discuss the estimators generated.

   For other estimators and more detailed investigations see \cite{Ve13,AO00};
for the DWR method see in particular \cite{BR03}.

From theoretical considerations it is clear that any good mesh for boundary or interior layers
must be anisotropic. Thus an adaptive procedure designed for problems with layers should include
an anisotropic refinement strategy. While several  anisotropic mesh
adaptation strategies do exist (see \cite{DF04} and its bibliography), all are more or less
heuristic. We do not know of  any strategy for convection-diffusion problems in two dimensions
where it is proved that, starting from some standard mesh, the refinement strategy is guaranteed
to lead to a mesh that allows robust error estimates.

Micheletti, Perotto and others \cite{DMP08,FMP04} combine SDFEM,\index{SDFEM} the DWR
method\index{DWR method} and anisotropic interpolation error estimates to get an \emph{a
posteriori} error estimate for some target functional. They then use this information to implement
a metric-based algorithm for mesh generation \cite{CLG07} that creates an ``optimal'' mesh. The
numerical results obtained are interesting but the second step of the approach has a heuristic
flavour.

Of course, strong results require an error estimator suitable for an anisotropic mesh. Several
authors use in the theoretical foundation of the estimators the {\it alignment measure}, introduced
by Kunert \cite {Ku99,Ku05}, see, for instance, \cite{Pi03,GSZ14,ZC14}. But the use of the
unknown, not computable alignment measure means that the initial mesh has more or less already
to reflect the anisotropy of the solution.

Kopteva designed different estimators ( residual \cite{Ko15,Ko17}, flux equilibration \cite{Ko18})
for anisotropic meshes, even in different norms (energy and maximum norm) without using any
alignment measure in proving upper error bounds.

Very important is the {\it robustness} of the estimators.
In \cite{San01} Sangalli proves the robustness of a certain  estimator
for the residual-free bubble method  applied to convection-diffusion
problems. The analysis uses the norm
\begin{equation}\label{eq3:3.533o}
||w||_{San}:=\|w\|_\varepsilon +\|b\cdot \nabla w\|_*,\quad\text{where}\quad \|\varphi\|_*=\sup
\frac{\langle\varphi,v\rangle}{\|v\|_\varepsilon}\,.
\end{equation}
Although Sangalli's approach is devoted to residual-free bubbles, the same analysis works for the
Galerkin method and the SDFEM. For the convection-diffusion problem, the residual error estimator
 is  robust with respect to the dual norm; see \cite{Ve04}.

Angermann's graph norm and the norm $||\cdot||_{San}$ above are defined only implicitly by an
infinite-dimensional variational problem and cannot be computed exactly in practice.
In \cite{San08}  Sangalli pointed out that the norm (\ref{eq3:3.533o}) seems to be  not optimal in
the convection-dominated regime. He proposes
an improved estimator that is robust with respect to his natural norm \cite{San05} for the
advection-diffusion operator,  but studied only the one-dimensional case.
 The relation to another new improved dual norm is in detail studied in \cite{DZ15}.

 Today the dual norm or its modification plays an important role in many papers on
robust a posteriori error estimation for convection-diffusion problems.

In \cite{TV15}, Tobiska and Verf{\"u}rth proved in the dual norm that the same robust a posteriori
error estimator can be used for a range of stabilized methods such as streamline diffusion, local
projection schemes, subgrid-scale techniques and continuous interior penalty methods. Nonconforming
methods are studied in \cite{ZC14}. Variants of discontinuous Galerkin methods are discussed in
\cite{ESV10,GSZ14,LPP16,SZ09,ZS11}. Vohralik \cite{Vo12} presents a very general concept of
a posteriori error estimation based on potential and flux reconstructions.

If one avoids to use a dual norm, often additional hypothesis are necessary to prove robustness.
See, for instance, \cite{JN13} for the analysis of an estimator for SDFEM in a natural SDFEM norm.

For reaction-diffusion problems
 the theoretical situation is clearer than for convection-diffusion
problems. The residual-based estimator  and the related estimator based on the
solution of auxiliary local problems are both robust with respect to the associated energy norm
\cite{Ve98b}. A modification of the equilibrated residual method
 of Ainsworth and Babuska is also robust for reaction-diffusion problems \cite{AB99}.

 For the flux reconstruction technique in the singularly perturbed case
equation \eqref{flux} corresponds to
\begin{equation}\label{flux2}
Res(v)=(f-cu_h+{\rm div}\,q,v)-(\varepsilon\nabla u_h-q,\nabla v)
\end{equation}
and one has carefully to estimate both terms of \eqref{flux2}. With some numerical
flux $q$ constructed, for instance, in \cite{CFPV09}, the first term yields the residual part
of the estimator
\[
\eta_{T,res}:=m_T\|f-cu_h+{\rm div}\,q\|_{0,T}
\]
with the same weight $m_T$ as in Verf{\"u}rth's estimator. The second term generates a
diffusive flux estimator $\eta_{DF}$, see \cite{CFPV09} for details. Flux equilibration
is also studied in \cite{AV11,AV14}, for nonconforming methods see \cite{ZCZ15}.

As in \cite{ZC14} is pointed out, the estimator of \cite{CFPV09} is not robust on
{\it anisotropic meshes}. A modification of that estimator is presented but the proof of
the upper bound uses the alignment measure. A recent result of Kopteva \cite{Ko18} for
flux equilibration on anisotropic meshes avoids the use of that ingredient.

Stevenson  proved in \cite{St05b} the uniform convergence of a special adaptive method for the
reaction-diffusion equation in the energy norm.

It is unclear that the energy norm is a
suitable norm for these problems because for small~$\varepsilon$ it is unable to distinguish
between the typical layer function of reaction-diffusion problems and zero. It would be desirable
to get robust \emph{a posteriori} error estimates in a stronger norm,  for instance, some
balanced norm or the $L_\infty$ norm.

The  first result with respect to the maximum norm is the \emph{a posteriori} error estimate of Kopteva
\cite{Ko07b} for the standard finite difference method on an arbitrary
rectangular mesh.  Next we sketch the ideas of \cite{DK16} for a posteriori error estimation for finite elements
of arbitrary order on {\it isotropic} meshes in the maximum norm.

Using the Green's function of the continuous operator with respect to a point $x$, the error in that point
can be represented by
\[
e(x)=\varepsilon^2(\nabla u_h,\nabla G)+(cu_h,G).
\]
For some $G_h\in V_h$ we obtain
\[
e(x)=\varepsilon^2(\nabla u_h,\nabla (G-G_h))+(cu_h,G-G_h).
\]
Integration by parts yields
\[
e(x)=\frac{1}{2}\sum_{T\in {\cal T}_h}\int_{\partial T}\varepsilon^2(G-G_h)n_T\cdot[\nabla u_h]
     +\sum_{T\in {\cal T}_h}(cu_h-f-\varepsilon^2\triangle u_h,G-G_h)_T.
\]
Choosing for $G_h$ the Scott-Zhang interpolant
 of $G$, one needs sharp estimates for $G$ to control
the interpolation error. These are
 collected in Theorem 1 of \cite{DK16}. Thus, one obtains finally with
 $l_{h}:=\ln(2+\tilde \varepsilon {\underline h}_{}^{-1} )$ (the constant $\tilde \varepsilon$ is
 of order $\varepsilon$ and $\underline h=\min h_T$)
 \begin{subequations}\label{Kop1}
\begin{align}
\|u-u_h\|_\infty &\le C\max_{T\in {\cal T}_h}(
   \min(\tilde\varepsilon,l_hh_T)\|[\nabla u_h]\|_{\infty,\partial T}\\
  &+ \min(1,l_hh_T^2\varepsilon^{-2})\|cu_h-f-\varepsilon^2\triangle u_h\|_{\infty, T})
.\nonumber
\end{align}
\end{subequations}

On {\it anisotropic meshes}, \index{anisotropic mesh} Kopteva also derives an a posteriori error estimator in the
maximum norm \cite{Ko15}, now for {\it linear} finite elements. Suppose that the
triangulation satisfies the maximum angle condition. Then the first result of \cite{Ko15}
gives
\begin{subequations}\label{Kop2}
\begin{align}
\|u-u_h\|_\infty &\le C\,l_h\max_{z\in \cal N}(
   \min(\varepsilon,h_z)\|[\nabla u_h]\|_{\infty,\partial \omega_z}\\
  &+ \min(1,h_z^2\varepsilon^{-2})\|cu_h-f\|_{\infty, \omega_z})
.\nonumber
\end{align}
\end{subequations}
Here $\omega_z$ is the patch of the elements surrounding some knot $z$ of the triangulation,
$h_z$ the diameter of  $\omega_z$. In a further estimator the second term of \eqref{Kop2},
which has isotropic character, is replaced by a sharper result with more anisotropic
nature.

To prove \eqref{Kop2} two difficulties arise. First, it is necessary to use scaled trace
bounds. Moreover, instead of using the Scott-Zhang interpolant of the Green's function (which
applicability is restricted on anisotropic meshes) Kopteva uses some standard Lagrange
interpolant for some continuous approximation of $G$. But the construction is based on the
following additional assumption on the mesh. Let us introduce $\Omega_1:=\{T: h_T\ge c_1\varepsilon\}$
and  $\Omega_2:=\{T: h_T\le c_2\varepsilon\}$ with some positive $c_1<c_2$. Then, the additional
assumption requires that the distance of $\Omega_1$ and $\Omega_2$ is at least some $c_3\varepsilon$
with $c_3>0$.

The last condition excludes an too abrupt change of the mesh size, typically for Shishkin meshes.
But other layer-adapted meshes satisfy that condition, for instance, Bakhvalov meshes or
Bakhvalov-Shishkin meshes.

\subsection{Time-dependent problems}
We start  from the problem
\begin{subequations}\label{eq3:4.1}
    \begin{align}
Lu := u_t -\varepsilon \Delta u + b\cdot\nabla u + cu
 &= f \quad\text{on }Q,  \label{eq3:4.1a} \\
%
 \text{ with initial-boundary}\text{
conditions }
 \nonumber\\
%
u(x,y,0) &=  s(x,y) \quad {\rm on\ } \Omega,
   \label{eq3:4.1b}\\
u(x,y,t) &=0 \quad {\rm on\ } \partial\Omega \times (0,T],
   \label{eq3:4.1c}
    \end{align}
\end{subequations}
assuming $c-1/2\nabla\cdot b\ge\beta>0.$

The weak form of \eqref{eq3:4.1} consists in finding $u\in L^2(0,T;H_0^1(\Omega))$ such that
$u_t\in L^2(0,T;H^{-1}(\Omega))$, $ u(\cdot,0)=s$ in $L^2(\Omega)$, and for almost every $t\in(0,T)$ and
$v\in H_0^1(\Omega)$
\begin{equation}\label{3.4.4.1}
 (u_t,v) +\varepsilon (\nabla u,\nabla v) + (b\cdot\nabla u + cu,v)
 = (f,v).
\end{equation}

In this subsection we mainly discuss two approaches for obtaining a posteriori error estimates: residual estimates
and estimates based on elliptic reconstruction.

Let us study the discretisation of \eqref{eq3:4.1} by linear finite elements in space and backward Euler in
time. Set $\tau_n:=t_n-t_{n-1}$, with every intermediate time $t_n$ we associate a shape-regular partition
${\cal T}_{h,n}$ and a corresponding finite element space $X_{h,n}$. Additionally, some transition condition
guarantees that the mesh at $t=t_n$ is not dramatically different from the mesh at $t=t_{n-1}$.\\
Then, the space-time discretisation reads: Find $u_h^n\in X_{h,n}$ such that
\begin{equation}\label{3.4.4.2}
   u_h^0=\pi_h s \quad  (L_2\,\, {\rm projection})
\end{equation}
and
\begin{subequations}
\begin{align}
(\frac{u_h^n-u_h^{n-1}}{\tau_n},v_h)&+\varepsilon (\nabla u_h^n,\nabla v_h) + (b^n\cdot\nabla u_h^n + c^nu_h^n,v_h)\nonumber\\
   &+\sum_{K\in {\cal T}_{h,n}}
   \delta_K\big(\frac{u_h^n-u_h^{n-1}}{\tau_n}-\varepsilon \Delta u_h^n + b^n\cdot\nabla u_h^n + c^nu_h^n,b^n\cdot\nabla v_h
   \big)_K\nonumber\\
   &=(f^n,v_h)+\sum_{K\in {\cal T}_{h,n}}
   \delta_K(f^n,b^n\cdot\nabla v_h)_K.
\end{align}
\end{subequations}
The choice $\delta_K=0$ yields standard Galerkin, otherwise we have SDFEM in space. Remark that we follow
\cite{Ve13}, where more general the $\theta$-scheme in time is studied ($\theta=1$ gives backward Euler).

The sequence $\{u_h^n\}$ defines a function $u_{h,\tau}$, piecewise affine on $(t_{n-1},t_n]$ with
$u_{h,\tau}(t_n)=u_h^n$. We equip the space
\[
X=\{u\in L^2(0,T;H_0^1(\Omega))\cap L^\infty(0,T;L^2(\Omega)):\,\, u_t+b\cdot\nabla u\in L^2(0,T;H^{-1}(\Omega))
\}
\]
with the norm defined by
\[
\|v\|^2_X:= {\rm sup}\|v(\cdot,t)\|_0^2+\int_0^T\|v(\cdot,t)\|^2_\varepsilon dt+
    \int_0^T\|( v_t+b\cdot\nabla v)(\cdot,t)\|_*^2 dt.
\]
As in the stationary case, residual based a posteriori error estimation uses the equivalence of residual and
error. More precisely, the norm of the residual $Res(u_{h,\tau)})\in L^2(0,T;H^{-1}(\Omega))$ defined by
\begin{equation}\label{3.4.4.3}
<Res(u_{h,\tau}),v>:=(f,v)-\big((\partial u_{h,\tau},v) +\varepsilon (\nabla u_{h,\tau},\nabla v) + (b\cdot\nabla u_{h,\tau} + cu_{h,\tau},v)\big)
\end{equation}
is equivalent to $\|u-u_{h,\tau}\|_X$, see Proposition 6.14 in \cite{Ve13}. Here we used
\[
\partial u_{h,\tau}:=\frac{u_h^n-u_h^{n-1}}{\tau_n} \quad {\rm in}\,\,(t_{n-1},t_{n}].
\]

Next we split the residual into a temporal and a spatial part. For simplicity, we restrict ourselves to time-independent
data, otherwise additional a data residual is necessary. Set
\[
<Res_h(u_{h,\tau}),v>:=(f,v)-\big((\frac{u_h^n-u_h^{n-1}}{\tau_n},v)+\varepsilon (\nabla u_{h}^n,\nabla v) + (b\cdot\nabla u_{h}^n + cu_{h}^n,v)
\big)
\]
and
\[
<Res_\tau(u_{h,\tau}),v>:=\varepsilon (\nabla (u_h^n-u_{h,\tau}),\nabla v) + (b\cdot\nabla(u_h^n- u_{h,\tau}) + (c\,(u_h^n-u_{h,\tau}),v).
\]
Consider first the spatial residual. \index{residual spatial} Analogously as in the stationary case, elementwise residuals $R_K$ and
edge residuals $R_E$ are introduced, and with the same weights $\alpha_K$ and $\alpha_E$ as in the stationary
case one generates the estimator
\[
\eta_h^n:=\big\{\sum_K\alpha_K^2\|R_K\|^2_{L^2(K)}+\sum_E\varepsilon^{-1/2}\alpha_E^2\|R_E\|^2_{L^2(E)}
\big\}^{1/2}.
\]
The direct estimation of the temporal residual yields the estimator \index{residual temporal}
\[
\hat\eta_h^n:=\{\|u_h^n-u_h^{n-1}\|^2_\varepsilon +\|b\cdot\nabla(u_h^n-u_h^{n-1})\|_*^2    \}^{1/2}.
\]
Together one gets a robust estimator, with, for instance, the following bound from above \cite{Ve13}:
\[
\|u-u_{h,\tau}\|_X\le C\{\|s-\pi_h s\|_0^2+\sum_n \tau_n((\eta_h^n)^2+(\hat\eta_h^n)^2)
\}^{1/2}.
\]
There exist also an estimate from below.

Unfortunately, $\|\cdot\|_*$ is not computable, and standard approaches, for instance inverse
inequalities, lead to estimates which are not robust with respect to $\varepsilon$. A similar but computable
estimator in \cite{AV14a} is also not fully robust.
 Therefore, Verf\"urth
introduces the auxiliary problem
\begin{equation}\label{3.4.4.4}
(\varepsilon (\nabla \tilde u_h^n,\nabla v_h) + \beta(\tilde u_h^n,v_h)=(b\cdot\nabla(u_h^n-u_h^{n-1}),v_h).
\end{equation}
One can show \cite{Ve13}, that finally one can replace
 $\|b\cdot\nabla(u_h^n-u_h^{n-1})\|_*$
by the computable quantity
\[
\|\tilde u_h^n\|_\varepsilon+\{\sum_K \alpha_K^2\| b\cdot\nabla(u_h^n-u_h^{n-1})+\varepsilon \Delta \tilde u_h^n-\beta \tilde u_h^n  \|^2_{0,K}
    +\sum_E\varepsilon^{-1/2}\alpha_E\|[n\cdot \nabla \tilde u_h^n]_E \|^2_{0,E}\}^{1/2}.
\]
The price for the final estimator is the need to solve the discrete stationary reaction-diffusion problem
\eqref{3.4.4.4} at each
time level.

Remark that in \cite{TV15} the authors study not only SDFEM but a wide range of stabilization methods, while
in \cite{AW17} for semilinear equations a fully adaptive Newton-Galerkin time stepping algorithm is designed.

Next we sketch the fundamental idea to use {\it elliptic reconstruction operators}\index{elliptic reconstruction} in a posteriori error
estimation for parabolic problems. Instead of comparing directly the exact solution with the numerical one, an
appropriate auxiliary function ${\cal R}u_{h,\tau}$ is defined. Then, we decompose the error into
\[
u-u_{h,\tau}=(u-{\cal R}u_{h,\tau})+({\cal R}u_{h,\tau}-u_{h,\tau}).
\]
The elliptic reconstruction ${\cal R}u_{h,\tau}$ is constructed in such a way that $u_{h,\tau}$ is the finite
element solution of an elliptic problem whose exact solution is ${\cal R}u_{h,\tau}$. Consequently,
${\cal R}u_{h,\tau}-u_{h,\tau}$ can be estimated by any available a posteriori error estimator for
elliptic problems. Moreover, $u-{\cal R}u_{h,\tau}$ satisfies a variant of the original PDE with a
right-hand side that can be controlled a posteriori. Then $u-{\cal R}u_{h,\tau}$ can be estimated using
well-known a priori estimates for the given time-dependent problem.

For simplicity let us start with the semi-discretisation of
\begin{subequations}\label{eq3:4.4.5}
    \begin{align}
u_t +Au
 &= f \quad\text{on }Q,   \\
 \text{ with initial-boundary}\text{
conditions }
 \nonumber\\
u(\cdot,0) &=  s \quad {\rm on\ } \Omega,\nonumber\\
u &=  0 \quad {\rm on\ } \partial\Omega \times (0,T]\nonumber
    \end{align}
\end{subequations}
by
\begin{equation}\label{eq3:4.4.6}
 (u_h)_t+A_h u_h=\pi_h f.
\end{equation}
Define the elliptic reconstruction of $u_h$ by \cite{MN03}
\begin{equation}\label{eq3:4.4.7}
   A({\cal R}u_h)=A_hu_h.
\end{equation}
It follows:
\[
   a({\cal R}u_h,v_h)=(A_hu_h,v_h)=a(u_h,v_h).
\]
That means: $u_h$ is the Ritz projection of ${\cal R}u_h$.

Next we derive an equation for $u-{\cal R}u_{h}$. Starting from
\[
(u-{\cal R}u_{h})_t+A(u-{\cal R}u_{h})=f-({\cal R}u_{h})_t-A_hu_h
\]
we get the error equation
\begin{equation}\label{eq3:4.4.8}
(u-{\cal R}u_{h})_t+A(u-{\cal R}u_{h})=f-\pi_h f+(u_h-{\cal R}u_{h})_t.
\end{equation}
Thus, we have the two properties required: $u_h$ is the Ritz projection of ${\cal R}u_{h}$,
and $u-{\cal R}u_{h}$ satisfies a variant of the original equation with a right-hand side
which can be controlled.

Based on available a posteriori information on $u_h-{\cal R}u_{h}$ in certain norms, from
\eqref{eq3:4.4.8} we can bound  $u-{\cal R}u_{h}$  using standard a priori estimates in
energy norms or in the $L^2$ norm or, alternatively, in the $L_\infty$ norm using estimates for the Green's
function of the given problem.

As an example for discretisation in space and time we follow  \cite{LM06} and consider backward
Euler in time:
\begin{equation}\label{eq3:4.4.9}
\frac{U^n-U^{n-1}}{\tau_n}+A^n U^n=\pi^n f^n.
\end{equation}
The continuous, piecewise linear interpolant of $\{U^N\}$ in time we denote by $U$, analogously
we define ${\cal R}U$ by the interpolate of  $\{{\cal R}^n\,U^n\}$. Here
\[
     A({\cal R}^n\,v)=A^n \,v.
\]
To derive the main parabolic error equation for $u-{\cal R}U$, we start from the discrete problem
in the form
\[
(\frac{U^n-U^{n-1}}{\tau_n},\pi^n\phi)+a(U^n,\pi^n\phi)-( \pi^nf^n,\pi^n\phi)=0
\]
or
\[
(\frac{U^n-\pi^nU^{n-1}}{\tau_n}+A^n\,U^n-\pi^nf^n,\pi^n\phi)=0.
\]
Because the quantity $(\frac{U^n-\pi^nU^{n-1}}{\tau_n}+A^n\,U^n-\pi^n f^n)$ lies in the finite element
space, it follows
\[
(\frac{U^n-\pi^nU^{n-1}}{\tau_n}+A^n\,U^n-f^n,\phi)=0 \quad \forall \phi\in H_0^1(\Omega).
\]
Equivalently
\[
(\frac{U^n-U^{n-1}}{\tau_n}+A^n\,U^n-\pi^nf^n,\phi)-(\frac{\pi^n U^{n-1}-U^{n-1}}{\tau_n},\phi)
=0 \quad \forall \phi\in H_0^1(\Omega).
\]
Introducing the elliptic projection we can write
\begin{equation}\label{eq3:4.4.10}
(\frac{U^n-U^{n-1}}{\tau_n},\phi)+a({\cal R}U^n,\phi)-(\pi^nf^n,\phi)-(\frac{\pi^n U^{n-1}-U^{n-1}}{\tau_n},\phi)
=0 .
\end{equation}
From that equation we substract the continuous problem
\[
  (u_t,\phi)+a(u,\phi)=(f,\phi)
\]
and obtain with the notation $\rho={\cal R}U-u$ and $\kappa={\cal R}U-U$ the error equation: for all
$\phi\in H_0^1(\Omega)$ it holds on $(t_{n-1},t_n)$
\begin{equation}\label{eq3:4.4.11}
(\rho_t,\phi)+a(\rho,\phi)=(\kappa_t,\phi)+a({\cal R}U-{\cal R}U^n,\phi)+(\pi^nf^n-f,\phi)+(\frac{\pi^n U^{n-1}-U^{n-1}}{\tau_n},\phi).
\end{equation}
So we can repeat the statement from above:\\
Based on available a posteriori information on $U-{\cal R}U$ in certain norms, from
\eqref{eq3:4.4.11} we can bound  $u-{\cal R}U$  using standard a priori estimates in
energy norms or in the $L^2$ norm, alternatively in the $L_\infty$ norm using estimates for the Green's
function of the given problem. See
 especially \cite{KL12,KL13,KL17}
for details concerning the technique described to obtain $L_\infty$ estimates for several discretizations
in time (backward Euler, Crank-Nicolson, dG). Backward Euler in combination with discontinuous Galerkin
in space is studied in \cite{CGM14}.

Remark that a different very general framework for robust a posteriori error estimation in unsteady problems
is presented in \cite{DEV13}. The authors use a special error measure (which cannot be computed easily in
practice), but obtain upper and lower bounds even in the case where the actual numerical scheme to
obtain $u_{h,\tau}$ need not to be specified. The upper bound for the error depends on an
convection-diffusion flux reconstruction and its local space-time equilibration, the lower bound
requires a local approximation property on this flux.

Goal-oriented a posteriori error control (the DWR method) is discussed in \cite{SB18}. For SDFEM in
space and discontinuous Galerkin in time the authors use the strategy first to dualize and then to
stabilize.

\section{Hints to further developments}
\subsection{Balanced norms}
Estimates in a balanced norm are popular since a work of Lin and Stynes 2012 \cite{LS12}. For a survey for second-order
reaction-diffusion problems, see \cite{Ro17}.
 Fourth-order problems are discussed in \cite{FR16}, higher order problems in \cite{FR20}.

 For time dependent problems the first approach to derive estimates in a balanced norm was published in \cite{BF22}.

\subsection{Flux-corrected transport}
 For the history of flux-corrected transport schemes see \cite{BJK18} and the references
given there. More recently,  the analysis of AFC schemes was presented  in three papers 2015-2017
by Barrenechea, John and Knobloch.

\subsection{Pressure-robust schemes}
Beginning with a paper of A. Linke 2014 (see \cite{Li14}) on pressure-robust schemes, a new aspect in the FEM analysis of incompressible flows, many papers were devoted on that topic
\cite{ALM17,LM16,LM16a,LMT16,JLM17,LLM17}. For a recent work see \cite{AK20}.

\bibliographystyle{plain}
\bibliography{ref,stynes,tobi,roo-g}

\end{document}